\documentclass[12pt]{article}
\usepackage{epsfig,latexsym,amsfonts,amsmath,amssymb}
\newtheorem{Theorem}{Theorem}
{Lemma}
{Proposition}
{Corollary}
\newtheorem{Remark}%[Theorem]
{Remark}

\def\tildeSobolev{\widetilde{H}^s(\Omega)}

\def\div{{\rm div}}

\begin{document}

\title %{\vspace{-10mm}
{On comparison of fractional Laplacians}

\author{Alexander I. Nazarov\footnote{St.Petersburg 
Department of Steklov Institute, Fontanka, 27, 
St.Petersburg, 191023, Russia
and St.Petersburg State University, 
Universitetskii pr. 28, St.Petersburg, 198504, Russia. 
E-mail: {al.il.nazarov@gmail.com}. Supported by RFBR grant 
20-01-00630.}
}

\date{}

\maketitle

\footnotesize

\noindent
{\bf Abstract.} For $s>-1$, $s\notin\mathbb N_0$, we compare two natural types of fractional Laplacians $(-\Delta)^s$, namely,
the restricted Dirichlet and the spectral Neumann ones. We show that for the quadratic form of their difference taken on the space 
$\tildeSobolev$ is
positive or negative depending on whether the integer part of $s$ is even or odd. For $s\in(0,1)$ and convex domains we prove also that the difference of these operators is positivity preserving on $\tildeSobolev$. This paper complements \cite{MN14} and \cite{MN16} where 
similar statements were proved for the spectral  Dirichlet and the 
restricted Dirichlet fractional Laplacians.

\normalsize

\bigskip\bigskip

\section{Introduction}

In recent decades a lot of efforts have been invested in 
studying nonlocal differential operators and nonlocal 
variational problems. Model operators here are various
fractional Laplacian (FLs for the brevity) $(-\Delta)^s$, mainly for $s\in(0,1)$.

Recall that the {\bf spectral Dirichlet and Neumann} 
FLs are the $s$th powers of 
conventional Dirichlet and Neumann Laplacian in the sense of 
spectral theory. In a Lipschitz bounded domain $\Omega$, 
they can be defined by corresponding quadratic forms 
$$
\aligned
Q_s^{\rm DSp}[u]\equiv &\, ((-\Delta_{\Omega})^s_{\rm DSp}u,
u):=\sum\limits_{j=1}^{\infty}\lambda_j^s|(u,\varphi_j)|^2; 
\\
Q_s^{\rm NSp}[u]\equiv &\, ((-\Delta_{\Omega})^s_{\rm NSp}u,
u):=\sum\limits_{j=0}^{\infty}\mu_j^s|(u,\psi_j)|^2, 
\endaligned
$$
where $\lambda_j$, $\varphi_j$ and $\mu_j$, $\psi_j$ are 
eigenvalues and (normalized) eigenfunctions of the 
Dirichlet and Neumann 
Laplacian in $\Omega$, respectively. Notice that $\mu_0=0$ 
and $\psi_0\equiv const$.

For $s\in(0,1)$ the 
domains of these quadratic forms are the classical
Sobolev--Slobodetskii spaces (see \cite[Ch. 4]{Tr} or 
\cite{DNPV})
\begin{equation}
\label{dom}
{\rm Dom}(Q_s^{\rm DSp})=\tildeSobolev;\qquad {\rm 
Dom}(Q_s^{\rm NSp})=H^s(\Omega)
\end{equation}
(we recall that
$$
\tildeSobolev=H^s(\Omega)\quad\mbox{if}\quad 0<s<1/2; \qquad 
\tildeSobolev\subsetneq H^s(\Omega)\quad\mbox{if}\quad 
s\ge1/2,
$$
see, e.g., \cite[4.3.2]{Tr}).

The first equality in (\ref{dom}) is proved in \cite[Lemma 
1]{MN14}; the 
proof of the second one is quite similar.

For $s>1$ the domains of spectral quadratic forms are more 
complicated but the following relations are always true:
\begin{equation*}
\tildeSobolev\subset{\rm Dom}(Q_s^{\rm DSp});\qquad 
\tildeSobolev\subset{\rm 
Dom}(Q_s^{\rm NSp}).
\end{equation*}

On the other hand, the quadratic form of {\bf restricted 
Dirichlet} FL is defined as follows:
$$
Q_s^{\rm DR}[u]\equiv ((-\Delta_{\Omega})^s_{\rm DR}u,
u):=\int\limits_{\mathbb R^n}|\xi|^{2s}|{\cal F}u(\xi)|^2\, 
d\xi
$$
where $\cal F$ is the 
Fourier transform
$$
{\cal 
F}{u}(\xi)=\frac{1}{(2\pi)^{\frac n2}}\int\limits_{\mathbb 
R^n} e^{-i\langle\xi,x\rangle}u(x)~\!dx.
$$
Corresponding domain is $\tildeSobolev$ for all $s>0$.

For $s\in(0,1)$ the following relation holds:
$$
Q_s^{\rm DR}[u]=c_{n,s}
\iint\limits_{\mathbb R^n\times\mathbb R^n} 
\frac{|u(x)-u(y)|^2}{|x-y|^{n+2s}}\, dx\,dy, $$
where
$$
c_{n,s}=2^{2s-1}\pi^{-n/2}\,\frac{\Gamma(\frac{n+2s}{2})}{
|\Gamma(-s)|}.
$$

\begin{Remark}
 Notice that for $s\in(0,1)$ the quadratic form of {\bf 
restricted Neumann} (or {\bf regional}) FL is
$$
Q_s^{\rm NR}[u]\equiv ((-\Delta_{\Omega})^s_{\rm NSp}u,
u):=c_{n,s}
\iint\limits_{\Omega\times\Omega} 
\frac{|u(x)-u(y)|^2}{|x-y|^{n+2s}}\, dx\,dy. 
$$
For some other types of fractional Laplacians see, e.g., 
\cite{MN19} and references therein.
\end{Remark}

The operators 
$(-\Delta_{\Omega})^s_{\rm DSp}$ and 
$(-\Delta_{\Omega})^s_{\rm DR}$ were compared in the sense 
of quadratic forms and in the pointwise sense
in \cite{MN14} ($s\in(0,1)$) and \cite{MN16} (for partial 
results see also \cite{ChSo}, \cite{F}, 
\cite{FG}, \cite{SV3}).

\begin{Theorem}[Theorem 2 in \cite{MN14} and Theorem 1 in 
\cite{MN16}]
\label{T:1}
Let $s>-1$ and $s\notin\mathbb N_0$. Suppose 
that\footnote{For $n=1$ and $s\le-\frac 12$ assume in 
addition that $(u,{\bf 1})=0$.} $u\in \tildeSobolev$, 
$u\not\equiv0$. Then the following relation 
holds:
\begin{equation*}
\aligned
Q_s^{\rm DSp}[u] > Q_s^{\rm DR}[u], & \quad {\rm if} & 
2k<s<2k+1,\ \ k\in\mathbb N_0;
%\label{pos_def}
\\
Q_s^{\rm DSp}[u] < Q_s^{\rm DR}[u], & \quad {\rm if} & 
2k-1<s<2k,\ \ k\in\mathbb N_0.
\endaligned
\end{equation*}
\end{Theorem}

\begin{Theorem}
\label{T:2}
\begin{enumerate}
 \item {\bf (Theorem 1 in \cite{MN14})}   Let $s\in(0,1)$, and 
let $u\in \tildeSobolev$, $u\ge0$, 
$u\not\equiv0$. Then the 
following relation 
holds in the sense of distributions:
\begin{equation*}
(-\Delta_{\Omega})^s_{\rm DSp}u > 
(-\Delta_{\Omega})^s_{\rm DR}u.
%\label{pos_pres>}
\end{equation*}

\item {\bf (Theorem 3 in \cite{MN16})}
Let $s\in(-1,0)$. Suppose 
that\footnotemark[1] $u\in \tildeSobolev$, $u\ge0$ in the 
sense of distributions, $u\not\equiv0$. Then the 
following relation holds:
\begin{equation*}
(-\Delta_{\Omega})^s_{\rm DSp}u < 
(-\Delta_{\Omega})^s_{\rm DR}u.
%\label{pos_pres<}
\end{equation*}
\end{enumerate}
\end{Theorem}

In this paper we prove similar results for the operators 
$(-\Delta_{\Omega})^s_{\rm DR}$ and 
$(-\Delta_{\Omega})^s_{\rm NSp}$. Since the domains of 
their quadratic forms are in general different, we 
consider them on the smaller domain $\tildeSobolev$.

\begin{Theorem}
\label{T:main1}
Let $s>-1$ and $s\notin\mathbb N_0$. Suppose 
that\footnote{For $s<0$ assume in 
addition that $(u,{\bf 1})=0$.} $u\in \tildeSobolev$, 
$u\not\equiv0$. Then the following relation 
holds:
\begin{eqnarray}
Q_s^{\rm DR}[u] > Q_s^{\rm NSp}[u], & \quad {\rm if} & 
s\in(2k,2k+1),\ \ k\in\mathbb N_0;
\label{pos_def1>}
\\
Q_s^{\rm DR}[u] < Q_s^{\rm NSp}[u], & \quad {\rm if} & 
s\in(2k-1,2k),\ \ k\in\mathbb N_0.
\label{pos_def1<}
\end{eqnarray}
\end{Theorem}

\begin{Remark}
Notice that a weaker inequality $Q_s^{\rm DSp}[u]\ge 
Q_s^{\rm NSp}[u]$ for $u\in \tildeSobolev$, $s\in(0,1)$, is a particular case of the 
well-known 
Heinz inequality \cite{H}. On the other hand, the inequality $Q_s^{\rm DR}[u]\ge 
Q_s^{\rm NR}[u]$ for $u\in \tildeSobolev$, $s\in(0,1)$, is trivial.
\end{Remark}

\begin{Theorem}
\label{T:main2}
 Suppose that $\Omega$ is convex. Let $s\in(0,1)$, and let 
$u\in \tildeSobolev$, $u\ge0$, $u\not\equiv0$. Then the 
following relation 
holds in the sense of distributions:
\begin{equation}
(-\Delta_{\Omega})^s_{\rm DR}u > 
(-\Delta_{\Omega})^s_{\rm NSp}u \qquad\mbox{in}\quad\Omega.
\label{pos_pres1}
\end{equation}
\end{Theorem}

The structure of our paper is as follows. In Section 2 we recall the basic facts on the 
generalized harmonic extensions related to fractional Laplacians of orders $\sigma\in(0,1)$ and $-\sigma\in(-1,0)$. Theorems \ref{T:main1} and \ref{T:main2} are proved in Section 3. Also we show that the assumption of convexity in Theorem \ref{T:main2} cannot be removed.

\section{Fractional Laplacians as D-to-N and\\ N-to-D operators}

It is a common knowledge nowaday that some of FLs of order $\sigma\in(0,1)$ are related to the so-called {\it harmonic extension in $n+2-2\sigma$ dimensions} 
and to the generalized Dirichlet-to-Neumann map (notice that for $\sigma=\frac 12$ it was known long ago). 

Let $u\in H^\sigma(\mathbb R^n)$ (in our consideration, we always assume that $u\in \widetilde H^\sigma(\Omega)$ is 
extended by zero to $\mathbb R^n$). 
In the pioneering paper \cite{CaSi}, it was shown that there exists a unique solution $w_\sigma^{\rm DR}(x,y)$ of the BVP in the half-space
$$-\div (y^{1-2\sigma}\nabla w)=0\quad \mbox{in}\quad \mathbb R^n\times\mathbb R_+;\qquad w\big|_{y=0}=u,
$$
with finite energy (weighted Dirichlet integral)
$$
{\cal E}_\sigma^{\rm 
DR}(w)=\int\limits_0^\infty\!\int\limits_{\mathbb{R}^n} 
y^{1-2\sigma}|\nabla w(x,y)|^2\,dxdy,
$$
and the relation
\begin{equation}
(-\Delta_{\Omega})^\sigma_{\rm 
DR}u(x)=-C_\sigma\cdot\lim\limits_{y\to0^+} y^{1-2\sigma }\partial_yw_\sigma^{\rm DR}(x,y)
\label{extension_DR}
\end{equation}
holds in the sense of distributions and pointwise at every 
point of smoothness of $u$. Here the constant $C_\sigma$ is given by
$$
C_\sigma:=\frac{4^\sigma\Gamma(1+\sigma)}{\Gamma(1-\sigma)}.
$$

Moreover, the function $w_\sigma^{\rm 
DR}(x,y)$ minimizes ${\cal E}_\sigma^{\rm DR}$ 
over the set 
$${\cal W}_\sigma^{\rm DR}(u)=\Big\{w(x,y)\,:\,
{\cal E}_\sigma^{\rm DR}(w)<\infty~,\ \ w\big|_{y=0}=u\Big\},
$$
and the following equality holds:
\begin{equation}
Q_\sigma^{\rm DR}[u]=\frac {C_\sigma}{2\sigma}\cdot {\cal E}_\sigma^{\rm 
DR}(w_\sigma^{\rm DR}).
\label{quad_DR}
\end{equation}

In \cite{ST} this approach was transferred to quite general situation. In particular, it was shown that for $u\in H^\sigma(\Omega)$ there is a unique solution $w_\sigma^{\rm NSp}(x,y)$ of the BVP in the half-cylinder
$$
-\div (y^{1-2\sigma}\nabla w)=0\quad \mbox{in}\quad 
\Omega\times\mathbb R_+;\qquad w\big|_{y=0}=u,
\qquad \partial_{\bf n}w\big|_{x\in\partial\Omega}=0
$$
(here $\bf n$ is the unit vector of exterior normal to $\partial\Omega$) having finite energy 
\begin{equation*}
{\cal E}_\sigma^{\rm 
NSp}(w)=\int\limits_0^\infty\!\int\limits_{\Omega} 
y^{1-2\sigma}|\nabla w(x,y)|^2\,dxdy,
\label{energy_NSp}
\end{equation*}
and the relation
\begin{equation}
(-\Delta_{\Omega})^\sigma_{\rm 
NSp}u(x)=-C_\sigma\cdot\lim\limits_{y\to0^+} y^{1-2\sigma }\partial_yw_\sigma^{\rm 
NSp}(x,y).
\label{extension_NSp}
\end{equation}
holds in the sense of distributions on $\Omega$ and pointwise at every 
point of smoothness of $u$.

Moreover, the function $w_\sigma^{\rm 
NSp}(x,y)$ minimizes ${\cal E}_\sigma^{\rm NSp}$ 
over the set 
$${\cal 
W}^{\rm NSp}_{\sigma,\Omega}(u)=\Big\{w(x,y)\,:\,
{\cal E}_\sigma^{\rm 
NSp}(w)<\infty~,\ \ w\big|_{y=0}=u\Big\},
$$
and the following equality holds:
\begin{equation}
Q_\sigma^{\rm 
NSp}[u]=\frac {C_\sigma}{2\sigma}\cdot {\cal E}_\sigma^{\rm 
NSp}(w_\sigma^{\rm 
NSp}).
\label{quad_NSp}
\end{equation}

In a similar way, one can connect FLs of order $-\sigma\in(-1,0)$ with the generalized Neumann-to-Dirichlet map. It was done in \cite{CDDS} for the spectral Dirichlet FL and in \cite{CbS} for the FL in $\mathbb R^n$ (and therefore for the restricted Dirichlet FL). Variational characterization of these operators was given in \cite{MN16}. We formulate this result for the operator\footnote{We emphasize that $(-\Delta_{\Omega})^{-\sigma}_{\rm DR}$ is not inverse to $(-\Delta_{\Omega})^\sigma_{\rm DR}$.
} $(-\Delta_{\Omega})^{-\sigma}_{\rm DR}$.

Let $u\in \widetilde H^{-\sigma}(\Omega)$ (for $n=1$ and $\sigma\ge\frac 12$ assume in 
addition that $(u,{\bf 1})=0$). We consider 
the problem of minimizing the functional 
$$
\widetilde{\cal E}_{-\sigma}^{\rm DR}(w)=
{\cal E}_\sigma ^{\rm DR}(w)\,-\,2\,\big(u,w\big|_{y=0}\big)
$$
over the set ${\cal W}_{-\sigma}^{\rm DR}$, that is closure of smooth functions on $\mathbb R^n\times\overline{\mathbb R}_+$ with bounded support, with respect to ${\cal E}_\sigma^{\rm DR}(\cdot)$. We notice that by the result of \cite{CaSi} the duality $\big(u,w\big|_{y=0}\big)$ is well defined.

If $n>2\sigma$ (this is a restriction only for $n=1$) then the minimizer is determined uniquely. Denote it by $w_{-\sigma}^{\rm DR}(x,y)$. Then formulae (\ref{extension_DR}) and (\ref{quad_DR}) imply the relations
\begin{equation}
Q_{-\sigma}^{\rm DR}[u]=-\frac {2\sigma}{C_\sigma}\cdot \widetilde{\cal E}_{-\sigma}^{\rm DR}(w_{-\sigma}^{\rm DR});\qquad
(-\Delta_{\Omega})^{-\sigma}_{\rm DR}u(x)=\frac {2\sigma}{C_\sigma}\,w_{-\sigma}^{\rm DR}(x,0)
\label{-D}
\end{equation}
(the second relation holds for a.a. $x\in\Omega$).

In case $n=1\le 2\sigma $ the minimizer $w_{-\sigma }^{\rm DR}(x,y)$
is defined up to an additive constant. However, by assumption 
$(u,{\bf 1})=0$ the functional $\widetilde{\cal E}_{-\sigma }^{\rm DR}(w_{-\sigma }^{\rm DR})$ does not depend on the choice of the constant, and the first relation in (\ref{-D}) holds. The second equality in (\ref{-D}) also holds if we choose the constant such that $w_{-\sigma }^{\rm DR}(x,0)\to0$ as $|x|\to\infty$.%\medskip

Notice that the function $w_{-\sigma}^{\rm DR}$ solves the Neumann problem in the half-space
\begin{equation*}
-\div (y^{1-2\sigma }\nabla w)=0\quad \mbox{in}\quad \mathbb R^n\times\mathbb R_+;\qquad \lim\limits_{y\to0^+} 
y^{1-2\sigma }\partial_yw=-u
%\label{eq:-CS}
\end{equation*}
(the boundary condition holds in the sense of distributions).
So, we obtain the ``dual'' Caffarelli--Silvestre 
characterization of $(-\Delta_{\Omega})^{-\sigma }_{\rm DR}$ as the Neumann-to-Dirichlet map.
\medskip

Now we introduce the ``dual'' Stinga--Torrea characterization of $(-\Delta_{\Omega})^{-\sigma }_{NSp}$ in almost the same way as it was done in \cite{MN16} for $(-\Delta_{\Omega})^{-\sigma }_{DSp}$. Namely, let $u\in H^{-\sigma}(\Omega)$ and let $(u,{\bf 1})=0$. Then the function $w_{-\sigma }^{NSp}(x,y)$ minimizing the functional
\begin{equation*}
\widetilde{\cal E}_{-\sigma}^{NSp}(w)={\cal E}_\sigma^{NSp}(w)\,-\,2\,\big(u,w\big|_{y=0}\big)
\end{equation*}
over the set 
$${\cal W}^{\rm NSp}_{-\sigma,\Omega}(u)=\Big\{w(x,y)\,:\,
{\cal E}_\sigma^{\rm 
NSp}(w)<\infty\Big\},
$$
is defined up to an additive constant. By assumption 
$(u,{\bf 1})=0$ the functional $\widetilde{\cal E}_{-\sigma }^{\rm NSp}(w_{-\sigma }^{\rm NSp})$ does not depend on the choice of the constant, and formulae (\ref{quad_NSp}) and (\ref{extension_NSp}) imply
\begin{equation}
Q_{-\sigma}^{NSp}[u]=-\frac {2\sigma }{C_\sigma }\cdot \widetilde{\cal E}_{-\sigma }^{NSp}(w_{-\sigma }^{NSp});\qquad
(-\Delta_{\Omega})^{-\sigma }_{NSp}u(x)=\frac {2\sigma }{C_\sigma }\,w_{-\sigma }^{NSp}(x,0)
\label{-N}
\end{equation}
(The second equality holds for a.a. $x\in\Omega$ if we choose the constant such that $w_{-\sigma}^{\rm NSp}(x,y)\to0$ as $y\to+\infty$).

Also the function $w_{-\sigma}^{\rm NSp}$ solves the Neumann problem in the half-cylinder
\begin{equation*}
-\div (y^{1-2\sigma }\nabla w)=0\ \ \mbox{in}\ \ \Omega\times\mathbb R_+;\quad \lim\limits_{y\to0^+} 
y^{1-2\sigma }\partial_yw=-u, \ \ \partial_{\bf n}w\big|_{x\in\partial\Omega}=0
%\label{eq:-ST}
\end{equation*}
(the boundary condition on the bottom holds in the sense of distributions).

\section{Proof of main results}

\noindent{\bf Proof of Theorem \ref{T:main1}}. We split the proof in three parts.

{\bf 1}. Let $s\in(0,1)$. 
For any $w\in{\cal W}_s^{\rm DR}(u)$ we have
$w\big|_{\Omega\times\mathbb R_+}\in{\cal W}^{\rm NSp}_{s,\Omega}(u)$. 
Therefore, relations (\ref{quad_DR}) and (\ref{quad_NSp}) provide
\begin{multline*}
Q_s^{\rm NSp}[u]=\frac {C_s}{2s}\cdot \inf\limits_{w\in{\cal 
W}^{\rm NSp}_{s,\Omega}(u)} {\cal E}_s^{\rm NSp}(w)
\le \frac {C_s}{2s}\, {\cal E}_s^{\rm NSp}(w_s^{\rm DR})\\
\le \frac {C_s}{2s}\, {\cal E}_s^{\rm DR}(w_s^{\rm DR})
=Q_s^{\rm DR}[u],
\end{multline*}
and (\ref{pos_def1>}) follows with the large sign.

Finally, the equality in (\ref{pos_def1>}) 
implies $\nabla w_s^{\rm DR}=0$ on $(\mathbb 
R^n\setminus\Omega)\times\mathbb R_+$. Since any 
$x$-derivative of $w_s^{\rm DR}$ solves the same equation  
in the whole half-space $\mathbb 
R^n\times\mathbb R_+$, it should be zero everywhere that is 
impossible for $u\not\equiv0$. 
\medskip

{\bf 2}. Let $s\in(-1,0)$. We define $\sigma=-s\in(0,1)$ and construct the extension $w_{-\sigma }^{\rm DR}$ as described in Section 2.

We again have
$w_{-\sigma }^{\rm DR}\big|_{\Omega\times\mathbb R_+}\in{\cal 
W}^{\rm NSp}_{-\sigma,\Omega}(u)$. 
Therefore, relations (\ref{-D}) and (\ref{-N}) provide
\begin{multline*}
-Q_s^{\rm NSp}[u]=\frac {2\sigma }{C_\sigma }\cdot \inf\limits_{w\in{\cal W}_{-\sigma,\Omega}^{\rm NSp}}\widetilde{\cal E}_{-\sigma}^{\rm NSp}(w)
\le \frac {2\sigma }{C_\sigma }\, \widetilde{\cal E}_{-\sigma}^{\rm NSp}(w_{-\sigma }^{\rm DR})\\
\le \frac {2\sigma }{C_\sigma }\, \widetilde{\cal E}_{-\sigma}^{\rm DR}(w_{-\sigma }^{\rm DR})=-Q_s^{\rm DR}[u],
\end{multline*}
and (\ref{pos_def1<}) follows with the large sign.
To complete the proof, we repeat the argument of the first part.\medskip

{\bf 3}. Now let $s>1$, $s\notin\mathbb N$. We put $k=\lfloor\frac {s+1}2\rfloor$ and define for $u\in\tildeSobolev$
$$
v=(-\Delta)^ku\in \widetilde H^{s-2k}(\Omega),\qquad s-2k\in (-1,0)\cup(0,1).
$$ Note that $v\not\equiv0$ if $u\not\equiv0$, and 
$$
(v,{\bf 1})={\cal F}v(0)=|\xi|^{2k}{\cal F}u(\xi)\big|_{\xi=0}=0.
$$
Then we have
$$
Q_s^{\rm DR}[u]=Q_{s-2k}^{\rm DR}[v],\qquad Q_s^{\rm NSp}[u]=Q_{s-2k}^{\rm NSp}[v],
%\endaligned
$$
and the conclusion follows from cases {\bf 1} and {\bf 2}.
\hfill$\square$\medskip

\noindent{\bf Proof of Theorem \ref{T:main2}}. We recall the representation formulae 
for $w_s^{\rm DR}$ and $w_s^{\rm NSp}$, see \cite{CaSi} and 
\cite{ST}, respectively:
$$
w_s^{\rm DR}(x,y)=const\cdot\int\limits_{\mathbb 
R^n}\frac {y^{2s}u(\xi)\,d\xi}{(|x-\xi|^2+y^2)^{\frac 
{n+2s}{2}}};
$$
$$
w_s^{\rm NSp}(x,y)=\sum\limits_{j=0}^{\infty} 
(u,\psi_j)_{L_2(\Omega)}\cdot{\cal 
Q}_s(y\sqrt{\mu_j})\psi_j(x),\qquad {\cal 
Q}_s(\tau)=\dfrac{2^{1-s}\tau^s}{\Gamma(s)}{\cal K}_s(\tau),
$$
where ${\cal K}_s(\tau)$ stands for the modified Bessel 
function of the second kind.

First of all, these formulae imply for $u\ge0$, 
$u\not\equiv0$
$$
\lim\limits_{y\to+\infty} w_s^{\rm DR}(x,y)=0;\qquad 
\lim\limits_{y\to+\infty} w_s^{\rm NSp}(x,y)=
(u,\psi_0)_{L_2(\Omega)}\cdot\psi_0(x)>0;
$$
the second relation follows from the asymptotic behavior 
(see, e.g., \cite[(3.7)]{ST})
$$
\aligned
{\cal K}_s(\tau)\sim &\, \Gamma(s)2^{s-1}\tau^{-s},\quad
\mbox{as}\quad \tau\to 0;\\
{\cal 
K}_s(\tau)\sim &\, \left(\dfrac{\pi}{2\tau}\right)^{\frac 
12}e^{ -\tau } \bigl(1+O(\tau^{-1})\bigr)\quad 
\mbox{as}\quad \tau\to+\infty.
\endaligned
$$

Next, for $x\in\partial\Omega$ we derive by convexity of 
$\Omega$
$$
\partial_{\bf n}w_s^{\rm DR}(x,y)= 
const\cdot\int\limits_{\mathbb 
R^n}\frac 
{y^{2s}\langle (\xi-x),{\bf n}\rangle 
u(\xi)\,d\xi}{(|x-\xi|^2+y^2)^{\frac 
{n+2s+2}{2}}}<0.
$$

Thus, the difference $W(x,y)=w_s^{\rm NSp}(x,y)-w_s^{\rm 
DR}(x,y)$ has the following 
properties in the half-cylinder $\Omega\times\mathbb R_+$:
$$
-\div (y^{1-2s}\nabla W)=0;\qquad W\big|_{y=0}=0;
\qquad W\big|_{y=\infty}>0;
\qquad \partial_{\bf n}W\big|_{x\in\partial\Omega}>0.
$$
By the strong maximum principle, $W>0$ in
$\Omega\times\mathbb R_+$. Finally, we apply the boundary 
point principle (the Hopf--Oleinik lemma, see \cite{KH}) to 
the function $W(x,t^{\frac 1{2s}})$ and obtain (cf. 
\cite[Theorem 1]{MN14})
$$
\liminf\limits_{y\to0^+} y^{1-2\sigma }\partial_yW(x,y)=\liminf\limits_{y\to 
0^+}\frac{W(x,y)}{y^{2s}}=\liminf\limits_{t\to 
0^+}\frac{W(x,t^{\frac 1{2s}})}{t}>0, \quad x\in\Omega.
$$
This completes the proof in view of 
(\ref{extension_DR}) and 
(\ref{extension_NSp}).\hfill$\square$%\medskip

\begin{Remark}
For non-convex domains the relation (\ref{pos_pres1}) does 
not hold in general. We provide corresponding 
counterexample.

Put temporarily $\Omega=\Omega_1\cup\Omega_2$ where 
$\Omega_1\cap\Omega_2=\emptyset$. If $u\ge0$ is a smooth 
function supported in $\Omega_1$ then easily
$(-\Delta_{\Omega})^s_{\rm NSp}u\equiv0$ in $\Omega_2$. On 
the other hand, $w_s^{\rm 
DR}(x,y)>0$ for all $x\in\mathbb R^n$, $y>0$, and the 
Hopf--Oleinik lemma gives $(-\Delta_{\Omega})^s_{\rm 
DR}u<0$ in $\Omega_2$.

Finally, if we join $\Omega_1$ with $\Omega_2$ by a small 
channel then the inequality $(-\Delta_{\Omega})^s_{\rm DR}u 
< (-\Delta_{\Omega})^s_{\rm NSp}u$ in 
$\Omega_2$ holds by 
continuity.

\end{Remark}

\end{document}